# On the history of nested intervals: from Archimedes to Cantor

G.I. Sinkevich

Saint Petersburg State University of Architecture and Civil Engineering Vtoraja Krasnoarmejskaja ul. 4, St. Petersburg, 190005, Russia

galina.sinkevich@gmail.com

*Abstract.* The idea of the principle of nested intervals or the concept of convergent sequences which is equivalent to this idea dates back to the ancient world. Archimedes calculated the unknown in excess and deficiency, approximating with two sets of values: ambient and nested values. G. Buridan came up with a concept of a point lying within a sequence of nested intervals. P. Fermat, D. Gregory, I. Newton, C. MacLaurin, C. Gauss, and J.-B. Fourier used to search for an unknown value with the help of approximation in excess and deficiency. In the 19th century, in works of B. Bolzano, A.-L. Cauchy, J.P.G. Lejeune Dirichlet, K. Weierstrass, and G. Cantor, this logical construction turned into the analysis argumentation method. The concept of a real number was elaborated in the 1870s in works of Ch. Méray, Weierstrass, H.E. Heine, Cantor, and R. Dedekind. Cantor's elaboration was based on the notion of a limiting point and principle of nested intervals. What we are going to consider now, is the genesis of this idea which dates back to the ancient world.

## Continuity of a set of real numbers

The property of a set of real numbers of being continuous, or complete, was stated in the form of several conceptions in the second half of the 19th century. The following properties of a set of real numbers were underlying each of these conceptions: 1. Each sequence of closed nested intervals has a nonempty intersection, and Archimedes' axiom is valid. 2. Each bounded subset has an upper bound. 3. Each Cauchy sequence converges, and Archimedes' axiom is valid. 4. Each infinite bounded subset has a limiting point (Bolzano-Weierstrass property). 5. Each bounded monotonic sequence converges.

Hilbert proved the equivalency of these properties at the turn of the 20th century. All the conceptions originated from ancient methods of proportions and exhaustion. However, more than 2,300 years had passed before the concepts of a number and continuity appeared. The history of this process is rich and interesting, and it keeps a lot of pages open. We are only going to unravel the first of the ideas here, i.e. the method of nested intervals or, which makes about as much sense, convergent sequences.

## Archimedes

The first theory of a real number was elaborated by Eudoxus, set forth by Euclid in his *Elements,* and consisted of two parts: the theory of proportions and the exhaustion method which was elaborated for incommensurable values and involved elaboration of a monotonic sequence of sums of known values approximating deficiency to the sought-for geometrical value. In his "correspondence with Dositheus[1]" cycle ('The Quadrature of the Parabola, 'On the Sphere and the Cylinder', 'On Conoids and Spheroids, 'On Spirals', and 'On the Measurement of a Circle'), in order to calculate the sought-for value, Archimedes created two sequences of values measured which approximates to the sought-for value in excess and deficiency.

---

[1] All these works were written in the form of letters to Dositheus of Pelusium, pupil of Conon of Samos.

In his letter to Eratosthenes (The Method of Mechanical Theorems), Archimedes noted that the exhaustion method was convenient when it was necessary to demonstrate the correctness of the foregone conclusion. In order to find the conclusion itself, Archimedes used the heuristic mechanical method of mathematical atomism. He presented intervals of lines consisting of material points, planar figures consisting of intervals, and bodies consisting of flats, and he determined distances between centers of gravity.

In his work 'On the Sphere and the Cylinder' Archimedes showed that the ratio of sums of two sequences was nearer and yet nearer to the unity. In the same work, Archimedes introduced the fifth assumption which was later named Archimedes axiom: "The larger of two unequal lines, surfaces or bodies is larger than the smaller one by a value which, if added to itself, can exceed any given value of those which may be in certain relation with one another" – it is translation from Russian [Archimedes, 1962, p. 97]. In Health translation: "Of unequal lines, unequal surfaces and unequal solids, the greater exceed the less by such a magnitude as, when added to itself, can be made to exceed any assigned magnitude among those which are comparable with [it and with] one another" [Archimedes, 1897, p.195]. In other words, the difference between partial sums of both sequences can be made arbitrary small.

In his work 'On Spirals', Archimedes stated the corollary as follows: "Since the area bounded by spirals is intermediate in magnitude between the circumscribed and inscribed figures, it follows that (1) *a figure can be circumscribed to the are such that it exceed the area by less than any assigned space*, (2) *a figure can be inscribed such that the area exceeds it by less than any assigned space*" [The works of Archimedes / Trans. T.L. Health, Cambridge Univ. Press, 1897, p. 368]. Thus, the sought-for value was located inside the converging range of measurement results. The length of the range was reducing by a value which was less than any predefined value, and the bounds of the range form sequences approximating to the sought-for value [Sinkevich, 2015a].

## 12th, 13th centuries

In universities of Paris of the 12th century and those of Oxford of the 13th century, lectures of professors involved reading classics, particularly Aristotle, and commenting upon them. As V.P. Zubov noted, the "mathematization" of Aristotle and "physicalization" of Euclid were characteristic of the 14th century [Zoubov, 1958, p. 622]. Scientists tried to separate the notions of a point, line, and surface from their physical interpretation (disputes of realists and nominalists). The continuum was considered as a whole geometric or physical object (the notion of a numerical continuum was formed but in the 19th century). The idea of a sequence and infinite sequence[2] appeared in works of Averroes and Albert the Great as a movement characteristic; and the difference between the kinematic and dynamic aspects was conceived. In the 14th century, this idea was developed in works of the Calculators, a group of scientists from Merton College in Oxford. Owing to them, the ancient tradition of assessment with the help of inequations was replaced by a new tradition of exact calculation, i.e. the equation. It was their merit that the notions of a "sequence", "intensity", and "instantaneous velocity" were introduced in the science, although not precisely defined. Richard Swineshead (also known as Suisset) was a representative of Merton school. His Book of Calculations was written somewhat in 1346. It was in this book, that he first introduced physical notions of a change and movement in mathematics, intension and remission of qualities (density and tenuity, force and resistance, quickness and slowness, warmth and coldness). Suisset introduced ordered scales of continuous changes with a correspondence between them [Shirokov, 1976, p. 134].

---
[2] res successive

## Buridan

Jean Buridan (abt.1300–1358), junior contemporary of Bradwardine, was William Ockhkam's student at Sorbonne. His concept of a continuum set forth in 'Questions on eight books of Aristotle' Physics' [Buridan, John, 1509, Subtilissimae Quaestiones super octo Physicorum libros Aristotelis, Paris. Rpr. 1964, as Kommentar zur Aristotelischen Physik, Frankfurt a. M.: Minerva] and in the treatise 'Quaestio de puncto' [Question on <the Nature of> Points] is of interest. Buridan emphasized the notion of a bound, noted the importance of the geometrical notion of 'contingence'. The link between this notion and the notion of 'contingence' in works of Lobachevsky was noted by V.S. Shirokov [Shirokov, 1978, p. 254]. Buridan like Bradwardine believed that the continuum consisted of points but its "point or instant were infinitely small"[3] [Buridan, 2006, p. 257]. "However, it does not follow from this that a point is a certain infinitely small value, because $\frac{1}{12}$ part of a continuum, which is not something actually infinitely small, can be also called a point" [Buridan, 2006, p. 305]. "Points in a continuum are mutually ordered; and the first point and consequently, the other ones may be separated. The infinity of intermediate points does not overrule the precedence due to which all points are assumed to actually exist one outside another and to be pretty ordered"[4] [Buridan, 2006, p. 311-313].

Buridan created a structure which served basis for an important elaboration in mathematics of the 19$^{th}$ century: a sequence of intervals each of which contained a continuum point: "if we take the first terminal point, we can identify the portion of the line which is nearest or nearer to it than all other portions which do not constitute part of the portion concerned and the portion concerned does not constitute part of other portions. However, if we take several such portions one of which does not constitute part of any other, then there will be no two portions equally near to the first point. Example: The first half, the first quarter, or the first eighth of the line is immediately adjacent to the first point. However, one of the quarters precedes all other quarters, and one of the eighths precedes all other eighths. Accordingly, we can say the same about the points, viz: whereas all points are located absolutely one outside the other one, one point must precede all other points"[5] [Buridan, 2006, p. 313].

---

[3] in infinitum parvum est punctus vel instans
[4] The issue of choice and ordering was considered by E. Zermelo in 1904.
[5] The covering idea was pitched by Bolzano in 1817, could be come across in Dirichlet's lectures of 1862, the covering lemma was stated by Heine in 1872, the covering theorem was proved by Borel in 1895 and extended by Lebesgue in 1898.

## 16th, 17th centuries

Since the Renaissance, the interest to Archimedes' tradition has been growing. Early in the 12th century, his work 'Measurement of the Circle' was known in Europe translated from Arabic into Latin, and in 1269, William of Moerbeke, student of Albert the Great, translated all major works of Archimedes from Ancient Greek. In 1544, works of Archimedes were published in Bazel in Greek and translated into Latin, which was conducive to public dissemination of his ideas. Scientists repeated the mechanics of his reasoning in both old and new problems. Archimedes' works had a great influence on Galilee and Stevin. One could come across the method of elaboration of sequences of ambient and nested values in interpolation of high order differences in 'Logarithmic Arithmetic' by G. Briggs in 1624 [Briggs, 1624], in P. Fermat's 'Méthodes de quadrature' [Fermat, 1891, p. 255-288]. P. Fermat told about them to a student of Galilee, B. Cavalieri; in interpolation formulas of D. Wallis, calculating number $\frac{4}{\pi}$ in 1656 [Wallis, 1656].

## Kepler. Cavalieri. Mathematical atomism, geometrical algebra as a heuristic method

Mathematical atomism re-emerged in the 16th century. I. Kepler (1571-1630) used indivisibles when cubing. In 1635, B. Cavalieri (1589-1647) published his 'Geometria indivisibilibus continuorum nova quadam ratione promota' (Geometry, developed by a new method through the indivisibles of the continua, 1635). Cavalieri highly appreciated the heuristic potential of the method of indivisibles. However, he rejected the possibility of algebraic reasoning of the method, which gave rise to contradictions noted by Galilee who was his teacher. J. Wallis, I. Barrow, and I. Newton used the productive, although contradictory, method of Cavalieri. MacLaurin wrote: "Cavalieri felt difficulties and advantages of the new method to the same extent. He spoke of it as if he had anticipated that the new science had to be put into an undisputable shape in order to satisfy the most fastidious geometricians" [MacLaurin, 1748/49, p. XLIX-L].

Thus, the two trends proposed in works of Archimedes have shaped: search for the desired result with the help of the method of indivisibles which had been developed in works of Kepler and Cavalieri, and reasoning of the so found result with the help of sequences converging in excess and deficiency, which was furthered by J. Gregory and C. MacLaurin. Please note that the goal of problem solving was to find a geometrical value (length, area, volume); in the 17th century, this process did not provide for a definition of a number as such.

## Gregory

*1668, James Gregory*, *'The True Areas of the Circle and the Hyperbola'*. A North Briton, James Gregory (1638-1675) was a predecessor of Newton in creation of analysis. In 1644-1667, he lived in Italy where Stefano degli Angeli, student of Cavalieri, studied in Padua. It was in the same place that he published his two works: 'The True Areas of the Circle and the Hyperbola' [Gregory, 1668a and 'General Sections of Geometry' [Gregorie, 1668b] where he applied the method of Archimedes to find curved surface areas, however, as he himself noted, in combination with a more convenient and brief method of indivisibles proposed by Cavalieri. Gregory's 'Geometry' already contained a proportion which enabled finding the length of a curve with the help of an element of the curve and main ideas of integral calculus. Gregory was the first to use the term 'convergence'.

In the work written in the same year, 'The True Areas of the Circle and the Hyperbola' [Gregory, 1668a], expressing all relations in proportions of inscribed and circumscribed figures, Gregory generated sequences approximating to the true value of the area of the hyperbolic segment in excess and deficiency. Thus, the tradition of Archimedes took another twist based on the method of indivisibles, which enabled a simplified work with proportional quantities (areas).

## Newton

The tradition of Archimedes underlies the calculation of quadratures in Newton's interpolation formulas. In 1669, Newton created an approximate method for algebraic equations using tangents [Newton, 1712]. In 1740, T. Simpson improved this method [Simpson, 1740]; in 1768, J. Mourraille (J.-R.-P. Mourraille, 1720-1808) noted that the convexity must be in place [A History of Algorithms, 1999, p. 179-183]; the convergence conditions of this method were described by J.B. Fourier in 1826[6] [Fourier, 1831]. This method involved the process of contraction of the interval which included the root of an equation. After S. Banach had set up the principle of contracting mappings in 1922 [Banach, 1922], in the middle of the 20$^{th}$ century, the tangent method was extended in works of L.V. Cantorovich based thereon (e.g. in [Cantorovich, 1949]).

In 1686, in 'Mathematical principles of natural philosophy' (Book I, Section I, 'Of the method of first and last ratio's of quantities') [Newton, 1729, p. 41-53], Newton stated 11 lemmas [the first-ever in mathematics] of the theory of limits including the conclusion from lemma IV as follows: "Hence if two quantities of any kind are any how divided into an equal number of parts; and those parts, when their numbers is augmented and their magnitude diminished *in infinitum*, have a given ratio one to the other; the first to the first, the second to the second, and so on in order; the whole quantities will be one to the other in that same given ratio" [ibidem, p.44, emphasis as in original]. Newton used this lemma to find the area of a curvilinear figure by way of comparing it with a known figure by coordinate. This was an algebraization of Cavalieri's method, however, Newton's indivisibles were replaced with variables tending to zero.

## MacLaurin

*Colin MacLaurin, 'Treatise of Fluxions', 1742.* The Scottish mathematician Colin MacLaurin (1698-1746) was a friend of Newton on whose recommendation he became professor of Edinburgh University, having in 1726 replaced James Gregory Jr. in his position[7]. After Newton died in 1727, his conception of analysis took the flak for the lack of clarity and reasoning. The crown of this criticism was 'The Analyst' by J. Berkeley published in 1734. MacLaurin decided to write his reasoning for Newton's Treatise of Fluxions. In 1742, his Treatise of Fluxions was published [MacLaurin, 1742]. It included a methodical and clear statement of Newton's method. This work was intended to be a textbook for youth. MacLaurin was trying to demonstrate the closeness of Newton's method and the ancient exhaustion method. An experienced educationalist, he found a perceivable image of convergence, having arranged convergent values – whether lengths, areas, or volumes – on a straight line: "we shall represent the circles and polygons by right lines, in the same manner as all magnitudes are expressed in the fifth book of the Elements" [ibid., p. 5]. Suppose the right lines *AB* and *AD* to represent the two areas of the circle that are compared together; and let *AP*, *AQ* represented any two similar polygons inscribed in these circles… When two variable quantities, *AP* and *AQ*, which always are in an invariable ratio to each other, approach at the same time to two determined quantities, *AB* and *AD*, so that they may differ less from them than by any assignable measure, the ratio of these limits *AB* and *AD* must be the same as the invariable ratio of the quantities *AP* and *AQ*: and this may be considered as the most simple and fundamental proportion in this doctrine, by which we are enabled to compare curvilineal spaces in some of the more simple cases."[8] [ibidem, p. 6].

---

[6] Published already after the death of Fourier in 1831.
[7] James Gregory (1666-1742), Professor of Mathematics and nephew of a prominent mathematician James Gregory (1738-1675), and brother of an astronomer and mathematician David Gregory (1659-1708).

MacLaurin sharpened this classical provision with the following elaboration: "In general, let any determined quantity *AB* be always a limit betwixt two variable quantities *AP*, *AQ*, which are supported to approach continually to it and to each other, so that the difference of either from it may become less than any assignable quantity, or so that the ratio of *AQ* to *AP* may become less than any assignable ratio of a greater magnitude to a lesser. Suppose also any other determined quantity *ab* to be always a limit betwixt the quantities *ap* and *aq*, and *aq* being always equal to *AQ* or less than it, let *ap* be either equal to *AP* or greater than it" [ibidem, p. 10]. MacLaurin demonstrated that limits *AB* and *ab* would equal each other[9]. In this event, all its values would be located on an interval, that is to say, would be demonstrable. In other words, if $AP < AB < AQ$, then $\frac{AP}{AQ}$ is limited to the relation of the larger value to the smaller, and if $AP \leq ap < ab < aq \leq AQ$, then $ab = AB$. P.P. Gaidenko [Gaidenko, 2011, p. 86] highly appreciated this MacLaurin's metaphor. MacLaurin was the first to introduce the term 'Archimedes' axiom'.

Thus, owing to this MacLaurin's metaphor, the method of convergent sequences gained an unambiguous demonstrativeness which in the following century enabled this method to raise to the next degree of generalization. His structure was still too far from the definition of a number – it was stated for values.

## 18th century

The tradition of the finite difference relations method which originated in works of B. Taylor, A. de Moivre, and I. Newton, and was furthered by J. Stirling, L. Euler, and J.L. Lagrange, continued in the 18th century as well. In his 'Elements of Geometry' of 1794, A.-M. Legendre [Legendre, 1823] defined each geometrical value as a number, and vice versa, found a geometrical value to correspond to each number. The elaboration of Ampère of 1806 proving the Lagrange mean value theorem is noteworthy. He developed continued inequalities with the estimated relation (mean fraction) in center and, decreasing the pitch, approximated to the sought-for value [Ampère, 1806]. There were no drawings or geometrical associations in Ampère's work.

---

[8] Reference to the following image:
https://ia601404.us.archive.org/BookReader/BookReaderImages.php?zip=/24/items/atreatisefluxio00conggoog/atreatisefluxio00conggoog_tif.zip&file=atreatisefluxio00conggoog_tif/atreatisefluxio00conggoog_0020.tif&scale=2&rotate=0

[9] Reference to the following image:
https://ia601404.us.archive.org/BookReader/BookReaderImages.php?zip=/24/items/atreatisefluxio00conggoog/atreatisefluxio00conggoog_tif.zip&file=atreatisefluxio00conggoog_tif/atreatisefluxio00conggoog_0021.tif&scale=2&rotate=0

The response error valuation problem which had been set up already by Galilee was developed in works of Lagrange and Laplace. In 1775/1776, Lagrange published his 'Memoir on applying the method of averaging out results of a large number of observations where advantages of this method for calculation of probabilities are considered and where various problems related to this issue are solved' [Lagrange, 1776]. In this memoir, Lagrange considered the probability of an error in the arithmetic average in various laws of allocation of errors.

In 1774, exploring the stability of the solar system, P.S. Laplace published his 'Memoir on the probabilities of reasons by event' [Laplace, 1774] where he assumed the density distribution of errors as a function: $\varphi(x) = \frac{m}{2} e^{-m|x|}$, $m > 0$. Gauss completed the creation of the theory of errors. In 1809, in his work entitled "Theory of the Motions of Celestial Bodies Moving about the Sun in Conic Sections' [Gauss, 1809, p. 258–260], Gauss stated the canonical orbit perturbations tracking theory (§177–182). He obtained the following result analyzing the distribution of observational errors based on this formula: "According to the elegant theorem first discovered by Laplace, integral $\int e^{-hh\Delta\Delta} d\Delta$ from $\Delta = -\infty$ to $\Delta = +\infty$ will equal $\frac{\sqrt{\pi}}{h}$ (where $\pi$ means the length of the semi-circumference of a unit radius) and our function will look as follows: $\varphi\Delta = \frac{h}{\sqrt{\pi}} e^{-hh\Delta\Delta}$. This function, which in reality cannot be found, strictly expresses the probability of an error: whereas potential errors are in any event confined within certain limits, the probability of errors outside these limits must equal zero, while our formula would always provide some value" [Gauss, 2013, p. 258–260]. Gauss reasoned this result minimizing the sum of a squared error. Thus, the principle of converging sequences was enriched by one more metaphor from the theory of errors. In Western reference materials, this Gauss' result is called the first statement of the squeeze theorem[10]. However, one more step had to be taken for the theorem on compressed variable to appear. They needed the notion of a converging sequence and continuous function for this purpose.

## Bolzano

*1817, Bernard Bolzano.* In 1817, in §§6 and 7 of his work entitled "Purely analytic proof of the theorem that between any two values which give results of opposite sign, there lies at least one real root of the equation" [Bolzano, 1817] (translated into Russian [Bolzano, 1955]), Bernard Bolzano introduced the convergence test of a sequence of partial sums. Having denoted the sum of the first $n, n+1, n+2, ..., n+r$ terms of series $F_1(x), F_2(x), F_3(x), ..., F_{n+r}(x), ...$, Bolzano stated and proved the theorem as follows: "If a series of values $F_1(x), F_2(x), F_3(x), ..., F_{n+r}(x), ...$ with a property that the difference between its $n^{th}$ term $F_n(x)$ and any subsequent $F_{n+r}(x)$ (no matter how remote it is from the first one) remains less than any prescribed value, and if $n$ is sufficiently large, then there is always a certain *constant value* (and only one value) to which the terms of this series are progressively approximating and which they can approach arbitrary close if we extend them sufficiently far" [Bolzano, 1817, p. 34-35]. I.e., following the idea of

---

[10] This is the theorem which in academic folklore is used to be called the theorem on two militiamen, carabineers, gendarmes, policemen, on a drunkard and two policemen, on three chords, squeezed theorem, pinched theorem, sandwich theorem. In textbooks on analysis of the 20$^{th}$ century, it is normally differentiated as an independent theorem.

Archimedes, Bolzano demonstrated that if the difference between partial sums (elements) of the sequence could be made arbitrarily small, then each such sequence would converge to a certain limit. Cauchy repeated this criterion without any proof in 1821. Since then, it has been called the Cauchy criterion, and sequences which meet this criterion were thereafter called Cauchy-Cantor sequences. This was the first contribution Bolzano made in the elaboration of the conception of a real number.

His second contribution made in the same work was the creation of a notion of the least upper bound $U$ of a linear numerical area (supremum). Bolzano considered certain property[11] $M$ fulfilled for $x<u$ at certain $u$. He elaborated a sequence of nested open intervals in the form of $\left(u+\frac{D}{2^m}+\frac{D}{2^{m+n}}, u+\frac{D}{2^m}+\frac{D}{2^{m+n-1}}\right)$, where property $M$ was not fulfilled so that $U-\delta<u<U+\varepsilon$. Thus, the upper bound $U$ must lie within the sequence of closed intervals in the form of $\left[u, u+\frac{D}{2^m}+\frac{D}{2^{m+n}}\right]$, although Bolzano did not point to this fact. He only noted that the upper bound need not necessarily belong to the set concerned (which means that it may not incorporate the largest element).

In such a manner, by intuition, Bolzano set the key lines for further development of the concept of a number: the principle of converging sequences, the principle of embedded sections, and the notion of supremum. In 1830s, Bolzano started creating the concept of a real number in terms of sections [Bolzano, 1931; Rykhlik, 1958]. Publishing of these manuscripts did not start until 1930. Unfortunately, Bolzano's involuntary dismissal from teaching and scientific isolation prevented his works from being disseminated. Although Bolzano had brilliantly mastered mathematical technique and was the first to introduce the analytical proof in analysis, his articles were fundamental and even philosophical in their nature rather than practical. Being in advance of their time, they gained popularity but half a century later. However, it should be noted that Cauchy, without referring to Bolzano, repeated both his definition of a continuous function [Sinkevich, 2012] and sequence convergence criterion with anticipatory reasoning of the geometrical progression and binomial expansion practically word for word.

The ingenious insight of Bolzano was appreciated only after his death owing to Herman Hankel who published and popularized his works [Sinkevich, 2014a]. Bolzano's works were highly appreciated by K. Weierstrass and G. Cantor.

## Cauchy

*1821 and 1823, Augustin Cauchy.* Mathematics in Napoleonic France was banished from universities and militarized. The emphasis was on training of military specialists and engineers. For this purpose, they established École Polytechnique and other schools of engineering. According to Riemann, in France, mathematics was computational, and in Germany, it was conceptual. The course of analysis given by Cauchy at École Polytechnique was brief and application-oriented. However, this course contained a systematic statement of the theory of limits, the most important analysis theorems devoted to continuous functions, differentiation, integration, and the theory of series. At the same time, Cauchy did not address the notion of irrational numbers considering them to be limits of sequences of rational numbers without defining the operations. The genius of Augustin Cauchy (1789-1857) was in the strict and clear generalization of achievements of his predecessors. Based on these achievements, he created a harmonious concept of analysis, which enabled him to obtain new results and form new sections in mathematics, e.g. the theory of residues.

---
[11] Bolzano meant the property of a function to maintain the sign.

In the Course of 1821, addendum III 'On computational solution of equations', Cauchy considered (without any references) the theorem to which Bolzano had devoted his work of 1817. Cauchy stated it as follows: "Let $f(x)$ be a real function of the variable $x$, which remains continuous with respect to this variable between the limits $x=x_0$ and $x=X$. If the two quantities $f(x_0)$ and $f(X)$ have opposite signs, we can satisfy the equation $f(x)=0$ with one or several real values of $x$ between $x_0$ and $X$" [Cauchy, 1821, p.378; Cauchy Engl. transl. 2009, p. 310]. Cauchy proved it having divided an interval into $m$ equal portions and selected such subinterval in this interval which has different signs on its extremes. Proceeding with this algorithm, Cauchy obtained a sequence of contracting intervals $\frac{X-x_0}{m^n}$ long, where the sequence of limiting points increases on the left and decreases on the right. Being different in sign, the values of the function in the limiting points approach zero. As the function is continuous, the common limit of sequences of the argument is the root of the equation. Bolzano proved this theorem in 1817 bisecting the interval and demonstrating with the help of this elaboration consistency of the existence of a supremum. As in many other cases [Sinkievich, 2013], Cauchy ingeniously simplified the idea of its proof and, in effect, formalized the squeeze theorem. Currently, the theorem on the existence of a root of a continuous function is known as Bolzano-Cauchy theorem which was first time ever stated for polynomials by Michele Rolle in 1690. Cauchy's reasoning was completely identical to that of Bolzano, including the preliminary resort to the geometrical progression. However, Cauchy introduced more convenient designations, and his statement was eloquent and concise. Stated by Cauchy, ideas that had been uttered by predecessor mathematicians took a rigid and orderly form, having developed into a comprehensive and relevant (for his day) course of analysis.

## Darboux

In 1875, Gaston Darboux in his 'Memoir on Discontinuous Functions' [Darboux, 1875], construct integral sums sequences converging to the integral right and left.

## Concepts of a number

Since 1860s, sweeping changes have taken place in mathematical analysis. The new needs of mathematics associated with the emergence of Fourier series, discontinuous functions, the need to categorize the discontinuing set, gave rise to a new reform in the analysis in the second half of the 19[th] century. The initiative was recaptured by German mathematicians led by Weierstrass. New concepts of a number, continuity and uniform continuity appeared; a covering lemma appeared as a consequence of the theorem on nested intervals. Introducing the new concept of a number through the notions of the upper bound and uniform convergence of series, Weierstrass developed the notion of continuity in the language of epsilon-delta; in 1869–1872, Charles Méray proposed the concept of a number based on the notion of converging sequences of Cauchy [Méray, 1869; Méray, 1872]; in 1872, Edward Heine proposed a concept based on fundamental (convergent) sequences of Cantor [Heine, 1872]; in 1872, Richard Dedekind proposed his concept based on the notion of section [Dedekind, 1872]; in 1872, Cantor wrote his first work in the theory of sets [Cantor, 1872] in which he introduced the notion of a limiting point. In this work, Cantor defined his fundamental sequences as follows: "When I am speaking of a numerical value in general, this happens particularly in the event that an infinite sequence of rational numbers $a_1, a_2, ..., a_n, ...$ is proposed set with the help of a certain law and possessing such property that difference $a_{n+m} - a_n$ becomes infinitely small with increasing $n$ no matter how large or small is the whole positive number $m$ or, in other words, that for an arbitrarily selected (positive rational) ε, there exists a whole number $n_1$ with a property that $|a_{n+m} - a_n| < \varepsilon$ if $n \geq n_1$ and $m$ is any positive whole number" [Cantor, 1872, p. 123-124]. "I would mention Book Ten of the 'Basics' by Euclid which remains the exemplar of the subject considered in these paragraphs" [Cantor, 1872, p. 127].

## Hilbert

The reform and arithmetization of mathematical analysis of the 19th century reached its completion with creation of concepts of a real number, continuity, and axiomatization of arithmetic [Sinkevich, 2014b]. Weierstrass asserted that a point on a straight line corresponded to each number. However, he was not sure whether the opposite was true. Cantor stated that there existed an unambiguous correspondence between numbers and points of a straight line but this could not be proved. With the help of the notion of a section, Dedekind demonstrated the continuity of a geometrical line and continuity of a set of real numbers [Sinkevich, 2015b]. The emergence of non-Euclidean geometry led to the need to analyze the axioms of geometry, the notion of continuity, and completeness. The axiomatic system of arithmetic appeared in works of Dedekind and Peano. The systems of axioms of arithmetic and geometry had to be generalized based on a converging evidence. In 1899, Hilbert introduced a new system of axioms having included the Archimedean and completeness axioms in it:

"*Measurement axiom or Archimedean axiom*. Let $A_1$ be any point upon a straight line between the arbitrarily chosen points *A* and *B*. Take the points $A_2, A_3, A_4, ...$ so that $A_1$ lies between *A* and $A_2$, $A_2$ between $A_1$ and $A_3$, $A_3$ between $A_2$ and $A_4$ etc. Moreover, let the segments $AA_1, A_1A_2, A_2A_3, A_3A_4, ...$ be equal to one another. Then, among this series of points, there always exists a certain point $A_n$ such that *B* lies between *A* and $A_n$.

*Completeness axiom*. To a system of points, straight lines, and planes, it is impossible to add other elements in such a manner that the system thus generalized shall form a new geometry obeying all of the five groups of axioms. In other words, the elements of geometry form a system which is not susceptible of extension, if we regard the five groups of axioms as valid" [Hilbert, 1950, p. 15].

Many attempts were made in the 19th century to build a geometry without the Archimedean axiom: in 1890, in his 'Foundations of Geometry', G. Veronese (1854–1917) proposed the concept of linear non-Archimedean continuum [Veronese, 1891]; works of O. Stolz [Stolz, 1881]. D. Hilbert [Hilbert, 1903] discussed this issue. According to Hilbert, "The completeness axiom gives us nothing directly concerning the existence of limiting points, or of the idea of convergence. Nevertheless, it enables us to demonstrate Bolzano's theorem by virtue of which, for all sets of points situated upon a straight line between two definite points of the same line, there exists necessarily a point of condensation, that is to say, a limiting point. From a theoretical point of view, the value of this axiom is that it leads indirectly to the introduction of limiting points, and, hence, renders it possible to establish a one-to-one correspondence between the points of a segment and the system of real numbers. However, in what is to follow, no use will be made of the "axiom of completeness" [Hilbert, (1899), 1923, p. 21; Hilbert, 1950, p.16].

# Kolmogorov

In the 20th century, investigations of Kolmogorov demonstrated that the completeness axiom can be replaced by the principle of embedded sections (Cauchy-Cantor fundamental sequences) together with the Archimedean axiom [Kolmogorov, 1946; Tikhomirov, 2014; Gladky, 2009; Rusakov, 2006]. In 1940s, Kolmogorov created an elaboration of real numbers as functions of a natural number [Kolmogorov, 1946].

Already late in the 19th century, new concepts of a number, continuity, and the theory of sets were included in the courses on the theory of functions of a real variable. In the 10th century, such course was given in Russia by S.O. Shatunovsky in Odessa, by his student G.M. Fichtengolz in St. Petersburg, and by N.N. Luzin, P.S. Aleksandrov, and A.N. Kolmogorov in Moscow. The original edition of the book of Aleksandrov and Kolmogorov entitled "Introduction into the theory of functions of a real variable" was published in 1933. Neither this edition nor the following two contained any axiomatic elaboration.

According to V.M. Tikhomirov, "In the autumn of 1954, A.N. Kolmogorov started giving the course of lectures entitled 'Analysis III' for the third-year students of the Department of Mechanics and Mathematics at which the author quoted herein studied. This was the first synthetic (i.e. incorporating several sections of mathematics) course in the history of the Department of Mechanics and Mathematics of Moscow State University. The program of the course was developed by A.N. Kolmogorov in 1940s and 1950s …

Kolmogorov introduced the axiomatic definition of real numbers as a totality which constitutes a complete linear ordered field. Having defined algebraic relations and the relation of order, he proceeded to the last axiom, i.e. completeness axiom. Kolmogorov called it the 'axiom of continuity'. He provided a number of axioms of continuity and proved their equivalency. These axioms were associated with names of those prominent mathematicians of the 19th century owing to whom analysis possessed its coherence, viz: Dedekind, Bolzano, Weierstrass, Cantor, and Cauchy.

These axioms were as follows:

## Section axiom

(A) Dedekind's axiom of section. *If set R is presented as a union of two non-vacuous non-overlapping sets X and Y where each element of X is less than any element of Y, then there exists element z with a property that $x \leq z \leq y$ for any $x \in X$ and $y \in Y$ (or, which means the same, either there is a maximum element in X or there is a minimum element in Y).*

## Upper bound axiom

(B) Bolzano's least upper (greatest lower) bound axiom. *Any set $S \subset R$ bounded from above (below) has the least upper (greatest lower) bound (i.e. element M (m) with such property that $x \leq M \ \forall x \in S \ (x \geq m \ \forall x \in S)$, and for any $\varepsilon > 0$ there exists element $\bar{\xi}(\varepsilon) \in S$ ($\underline{\xi}(\varepsilon) \in S$) with such property that $\bar{\xi}(\varepsilon) > M - \varepsilon$ $(\underline{\xi}(\varepsilon) < m + \varepsilon)$.*

## Limiting point axiom

(C) Weierstrass' limiting point axiom. *Any bounded sequence $\{x_k\}_{k \in N}$ of elements from R has a limiting point (i.e. element $\xi \in R$, and any ε-neighborhood of this element contains an element of a sequence other than ξ).*

## Converging subsequence axiom

(D) Weierstrass's converging subsequence axiom. *A converging subsequence can be selected from any bounded sequence of elements from R.*

## Monotonic sequence axiom

(E) Bolzano's monotonically increasing (decreasing) sequence. *A bounded monotonically increasing (decreasing) sequence of elements from R has a limit.*

## Nested intervals axiom

As a consequence, Cantor's axiom of nested intervals comes: a sequence of nested intervals $\Delta_n = [a_n, b_n] \subset R$, $\Delta_1 \supset \Delta_2 \supset ...$ with lengths tending to zero (i.e. $b_n - a_n \to 0$) has the only common point ξ (i.e. $\xi \in \Delta_n \ \forall n$)" [Tikhomirov, 2014, p. 151–152].

Dedekind stated his axiom in 1872 [Dedekind, 1872], although the prescience of the notion of a section was stated in works of Bolzano in 1830s [Bolzano's Schrifte, 1931, Sinkiewich, 2013]. Bolzano stated his axiom on the existence of the least upper bound in 1817 [Bolzano, 1817]; the lack of the theory of a real number enabled Bolzano to only prove the consistency of the assumption that the upper bound exists. Cantor introduced the notion of a limiting point in 1872 [Cantor, 1872], and it was developed in lectures of Weierstrass [Weierstrass, 1989], also in Russian [Sinkevich, 2014c]. The notion of a converging subsequence was first stated by H.E. Heine in 1872 based on Cantor's idea and talks with Weierstrass [Heine, 1872]. As a conclusion from the axiom on embedded sections, the covering lemma comes.

After the work of Dedekind had been published, R. Lipschitz wrote to him: "I cannot deny relevancy of your definition; I just think that it differs only in the form of expression from what the ancients ascertained, not in the essence. The only thing I can say, is that I believe the definition provided by Euclid (Book V, definition 4) and your definition are equally satisfactory". Responding to Lipschitz, Dedekind, however, insisted that "Euclidian principles alone without resort to the continuity principle which is not contained in them cannot substantiate the perfect theory of real numbers as ratios of values ... And vice versa, owing to my theory of irrational numbers, an epitome of a continuous area has been created, which exactly on this reason can characterize any relation of values by certain numerical individuum contained therein" [Gaidenko, 2011, p. 82].

According to Dedekind, the mankind gradually ascended up the stairway of meanings (Treppenverstand), thoroughly partitioning the array of thoughts on which the rules of numbers rest [Dedekind, 1888, p. 3]. Considering the history of the method of nested intervals, we can see the ascent up the stairway of meanings from the ancient world to nowadays. The method of indivisibles existed for the sake of search; the exhaustion method, for Archimedes to assert; the notion of a sequence was conceived by the Calculators and took logical form in works of scholastics. The Renaissance aroused Kepler's and Cavalieri's interest to indivisibles again and encouraged Gregory to synthetize methods and generate the notion of convergence. Newton and Leibniz invented mathematics of variables, i.e. Calculus. MacLaurin was the first to emphasize the role of Archimedes axiom. He also created a metaphor, a demonstrable image of sequences converging to one another, as a sequence of nested intervals. The theory of measurement errors set up by Galilee, furthered in works of Laplace, and completed by Gauss, produced another metaphor: distribution of errors under a normal law. The development of the notion of a function and particularly continuous function in works of Bolzano determined the lines of development of the concept of continuity with the help of the notions of the upper bound, converging sequences, and section. Mathematical analysis of the early $19^{th}$ century was based on Cauchy's theory of limits sufficient at that time and his theorem on continuous functions. Based thereon, a couple of concepts of a number and continuity appeared in the second half of the $19^{th}$ century: those of Méray, Weierstrass, Cantor, and Dedekind. Hilbert found a solution to the differences in these concepts, having proved their equivalency; and in the $20^{th}$ century, Kolmogorov developed a unifying concept of a real number.

I would conclude this article with a quote from Dedekind: "The greatest and most fruitful success in mathematics and in other sciences is often achieved owing to new notions invented and introduced at the time when we are coerced to do so often addressing complex phenomena which can be explained with the help of first notions in a very complex form" [Dedekind, 1888, p. 5].

## References

A History of Algorithms. From the Pebble to the Microchip. 1999 /ed. Chabert L.-L. Springer.

Ampère, A., 1806. Recherches sur quelques points de la théorie des fonctions dérivées qui conduisent à une nouvelle démonstration de la série de Taylor et à l'expression finie des termes qu'on néglige lorsqu'on arrete cette série à un terme quelconque/Mémoir par M. Ampère, Répétiteur à l'Ecole Politechnique // Journal de l'Ecole Politechnique. Cahier 13, 148–181.

Archimedes, 1897. The works of Archimedes / Trans. T.L. Health, Cambridge Univ. Press. , 1897, p.195.

Archimed, 1962, in Russian: Sochinenija / Trans., preface and comments by I.N. Veselovski. Moscow: GIFML.

Banach S. 1922. Sur les operations dans les ensembles abstraits et leur application aux équations intégrales //Fundamenta Mathematicae, 3, 133–182.

Bolzano, B. 1817. Rein analytischer Beweis des Lehrsatzes, daß zwischen zwey Werthen, die ein entgegengesetztes Resultat gewähren, wenigstens eine reelle Wurzel der Gleichung liege. – Prag: Gottlieb Haase.

Bolzano, B., 1931. Zahlentheorie / Bernard Bolzano's Schriften, Band. 2. Praha, 1931.

Bolzano, B., 2004. The Mathematical Works of Bernard Bolzano, ed. and tr. S. B. Russ. Oxford University Press, Oxford.

Briggs H. 1624. Arithmetica logarithmica. Londini, 1624.

Buridan, J, 2006 (1350). Quaestio de puncto / translated into Rusian by V.P. Zubov: Traktat 'O Tochke'// V. Zubov. Iz istorii mirovoy nauki. Izbrannyie trudy 1921-1963. S.-Petersburg: Aleteja, 311-347.

Cantor G. 1872. Über die Ausdehnung eines Satzes aus der Theorie der trigonometrischen Reihen // Math. Ann., Bd.5, 123-132.

Cauchy, A.-L. 1821. Course d'Analyse de l'Ecole Royale Politechnique. Analyse Algébrique // Oeuvres. Ser. 2, t. 3.

Cauchy's Cours d'analyse, 2009 / An annotated English translation by R.E. Bradley, C.E. Sandifer. Springer, 2009.

Darboux, G., 1875. Mémoire sur les fonctions discontinues //Annales de l'École Normale. 1875. – 2-e Série. – Tome IV. – p. 57–112.

Dedekind, R. 1872. Stetigkeit und irrationale Zahlen. Braunschweig: Vieweg.

Dedekind, R. 1888. Was sind und was sollen die Zahlen? 1 Auflage, Braunschweig: Vieweg.

Fermat, P.,1891. Oeuvres. Paris, t.1.

Fourier, J.B.J., 1831. Analyse des équations determinées. Première partie. Paris: Chez Firmin Didot frères, libraries, Rue Jacob 24.

Gaidenko, P.P., 2011. Stanovlenije novoevropeiskogo estestvoznanija: preodolenije paradoxov aktualno-beskonechnogo (Formation of modern European science: overcoming the paradoxes of actual infinity) // Metafizika, No 1, 65-87.

Gauss, C.F., 2013 (1809).Theory of the motion of the heavenly bodies moving about the sun in conic sections. With an appendix /Translated in English by C. H. Davis. Boston. – 1857. – 268 p. India: Pranava Books. Reprint 2013. – 408 p. – P. 258–260., c. 258–260.

Gladkij, A.V., Koziorov, J.N., 2009. Dejstvitelnyje chisla kak posledovatelnosti obyknovennych drobej (Teorija dejstvitelnych chisel po Kolmogorovu). (Real numbers as a series of fractions (real number theory as on Kolmogorov)) // Matematika v vyshem obrasovanii, 7, 21-38.

Gregory J. 1668a. Vera circuli et hyperbolae quadratura cui accedit geometria pars vniuersalis inseruiens quantitatum curuarum transmutationi & mensurae. Authore Iacobo Gregorio Abredonensi. Padua, 1668.

Gregorie J. 1668b. The Universal Part of Geometry devoted to the transmutation and measurement of curved quantities. English trans. by Andrew Leahy. http://math.knox.edu/aleahy/gregory/WORKING/gpu.html

Heine, E. 1872. Die Elemente der Functionenlehre // J. reine angew. Math., 74, 172–188.

Hilbert, D., 1903. Grundlagen der geometrie. Zweite Auflage. Leipzig.

Hilbert, D., 1923. Osnovanija geometrii (the foundations of Geometry) / Russian trans. By A.V. Vasiljev. Petrograd.

Hilbert, D., 1950, (1902), The Foundations of Geometry [Grundlagen der Geometrie] / English translation by E.J. Townsend (2nd ed.), La Salle, IL: Open Court Publishing.

Kantorovich, L.V., 1949. O metode Newtona (On the Newton's method) // Sbornik rabot po priblizchennomu analizu Leningradskogo otdelenija Instituta, Trudy MIAN USSR, Moscow-Leningrad, 28, 104-144.

Kolmogorov, A.N., 1946. K obosnovaniju teorii deistvitelnych chisel (Substantiation of the theory of real numbers)//Uspechi matematicheskich nauk, 1, 217-219.

Lagrange, J.L., 1776. Mémoire sur l'utilité de la méthode de prendre le milieu entre les résultats de plusiers observations, dans lequel on examine les avantages de cette méthode par le calcul des probabilités, et où l'on résout différents problèmes relatifs à cette matière //Misc. Taur., pour 1770-1773, v.5 (1776), p. 167-232 ; Œuvres, t. 2, 173-234.

Laplace, P., 1774. Mémoire sur la probabilité des Causes par les Événements // Mémoires l'Academie Royale des Sciences. Paris, 6, 621–656.

Legender, A. M., 1823. Les Éléments de géométrie. Par A. M. Legender, Avec Additions Et Modifications. – Paris.

MacLaurin, C. A., 1742. Treatise of Fluxions in two books by Colin MacLaurin, A.M., Professor of Mathematics in the Univesity of Edinburg and Fellow of the Royal Society. Edinburg: Printed by T.W. and T. Ruddmans.

Maclaurin, C. A., 1748/49.Traité des fluxions, trad. par Pezenas.

Méray, Ch.,1869. Remarques sur la nature des quantités définies par la condition de servir de limites à des variables données // Revue des Sociétés savantes, Sci. Math. phys. nat. (2) 4, 280–289.

Méray, Ch., 1872. Nouveau précis d'analyse infinitesimale. Publication: F. Savy. XXIII

Newton, I., 1712. Analysis per aequationes numero terminorum infinitas // Commercium epistolicum D. Johannis Collins et aliorum de analysi promota. Sent by Dr. Barrow to Mr. Collins in a Letter dated July 31.1669. Edita Londini, 1712. - P. 3–20.

Newton, I., 1729 (1687) The Mathematical Principles of Natural Philosophy / Translated into English by A. Motte. London. v. I.

Rusakov, A.A., Chubarikov V.N., 2006. O dvuch podchiodach k obosnovaniju vezhchestvennych chisel (Two approaches to the justification of real numbers) // Matematika v vysshem obrasovanii, 4, 37-44.

Rychlik, K., 1958. Teorija vezhshestvennych chisel v rukopisnom nasledii Bolzano (The theory of the real numbers in the manuscript heritage of Bolzano) Istoriko-matematicheskije issledovanija, XI, 515–532.

Simpson, Th., 1740. Essay on several subjects in speculative and mixed mathematics. London.

Sinkevich, G., 2012. K istorii epsilontiki (On the history of epsilontics) // Matematika v vyshem obrazovanii, 10, 149–166. + Sinkevich G. On the history of epsilontics// http://arxiv.org/abs/1502.06942

Sinkiewich, G. I., 2013. Historia dwóch twierdzeń analizy matematycznej: M. Rolle, B. Bolzano, A. Cauchy // Dzieje matematyki polskiei II. Praca zbiorowa pod redakcją Witolda Więsława. Institut Matematyczny Uniwersytetu Wrocławskiego. Wrocław, 165 – 181. + Sinkevich G. Rolle' Theorem and Bolzano-Cauchy' Theorem from the end of the 17th century to K. Weierstrass epoch. http://arxiv.org/abs/1503.03118

Sinkevich G. 2014a. Rasprostranenie i vlijanije idei Bolzano na razvitije analisa XIX veka (The spread and influence of ideas of Bolzano on the analysis development of the XIX century). The Abstracts of International Conference "Infinite dimensional analysis, stochastics, mathematical modeling: new problems and methods" in Peoples' Friendship University of Russia (Moscow) on December 15-18, 2014, 436–438.

Sinkevich G. 2014b. Concepts of a Numbers of C. Méray, E.Heine, G. Cantor, R. Dedekind and K. Weierstrass // Technical Transactions. Kraków, 1-NP, 211–223.

Sinkevich, G. 2014c. The Notion of Connectedness in Mathematical Analysis of XIX century / G. Sinkevich // Technical Transactions. Kraków. 1-NP, 195-209.

Sinkevich, G., 2015a. Archimedes: pis'ma k Dositheus i axioma polnoty (Archimedes: correspondence with Dositheus and Axiom of Completeness)// Proceedings of International Conference "Infinite dimensional analysis, stochastics, mathematical modeling: new problems and methods" in Peoples' Friendship University of Russia (Moscow) on December 15-18, 2014, 366-370.

Sinkevich, G., 2015b. On the History of Number Line. http://arxiv.org/abs/1503.03117

Shirokov, V.S., 1976. O 'Knige vychislenij' Richarda Suisseta (On the 'Book of Calculations' by Richard Suisset) // Istoriko-matematicheskije issledovanija, XXI, 129–142.

Shirokov, V.S., 1978. Infinitesimalnayja koncepcija Buridana (Buridan's conception of infinitesimal)// Istoriko-matematicheskije issledovanija, XXIII, 250–269.

Stolz O. 1881. B. Bolzano's Bedeutung in der Geschichte der Infinitesimalrechnung. – Math. Ann. Bd. 18, 255-279.

Tichomirov, V.M. 2014. Axiomaticheskij metod I teorija dejstvitelnych chisel v lektsijach A.N. Kolmogorova (Axiomatic method and the theory of the real numbers in AN Kolmogorov's lectures) // Matematika v vysshem obrasovanii, 12, 149-154.

Veronese, G. 1891. Fondamenti di geometria a più dimensioni e a più specie di unità rettilinee esposti in forma elementare. Padova: Tipografia del Seminario.

Wallis. 1656. Arithmetica infinitorum, Oxoniae.

Weierstrass, K., 1989. Ausgewählte Kapitel aus der Funktionenlehre. Vorlesung gehalten in Berlin 1886 mit der Akademischen Antrittsrede, Berlin 1857 und drei weiteren Originalarbeiten von K. Weierstrass aus den Jahren 1870 bis 1880/86. Teubner-Archiv für mathematic. Band 9. Reprint 1989.

Zubov, V. P., 1958. Traktat Nicolas Oresme 'O konfiguratsii kachstv' (Nicolas Oresme's Tractatus de configurationibus qualitatum et motuum// Istoriko-matematicheskije issledovanija. Moskwa : Nauka, XI, 601–635.

Zoubov, V. (ed.), 1961, "Jean Buridan et les concepts du point au quatorzième siècle," Medieval and Renaissance Studies, 5: 63–95.